\pdfoutput=1
\RequirePackage{ifpdf}
\ifpdf 
\documentclass[pdftex]{sigma}
\else
\documentclass{sigma}
\fi

\begin{document}

\allowdisplaybreaks

\renewcommand{\PaperNumber}{004}

\FirstPageHeading

\ShortArticleName{On a Lie Algebraic Characterization of Vector Bundles}

\ArticleName{On a Lie Algebraic Characterization of Vector Bundles}

\Author{Pierre B.A.~LECOMTE, Thomas LEUTHER and Elie ZIHINDULA MUSHENGEZI}

\AuthorNameForHeading{P.B.A.~Lecomte, T.~Leuther and E.~Zihindula Mushengezi}

\Address{Institute of Mathematics, Grande Traverse 12, B-4000 Li\`ege, Belgium}
\Email{\href{mailto:plecomte@ulg.ac.be}{plecomte@ulg.ac.be}, \href{mailto:thomas.leuther@ulg.ac.be}{thomas.leuther@ulg.ac.be}, \href{mailto:eliezihindula@ulg.ac.be}{eliezihindula@ulg.ac.be}}

\ArticleDates{Received September 23, 2011, in f\/inal form January 23, 2012; Published online January 26, 2012}

\Abstract{We prove that a vector bundle $\pi : E \to M$ is characterized by the Lie algebra generated by all dif\/ferential operators on $E$ which are eigenvectors of the Lie derivative in the direction of the Euler vector f\/ield. Our result is of Pursell--Shanks type but it is remarkable in the sense that it is the whole f\/ibration that is characterized here. The proof relies on a~theorem of [Lecomte~P.,
\textit{J.~Math. Pures Appl.~(9)} \textbf{60} (1981), 229--239] and inherits the same hypotheses. In particular, our characterization holds only for vector bundles of rank greater than~$1$.}


\Classification{13N10; 16S32; 17B65; 17B63}

\vspace{-2mm}

\section{Introduction}

Algebraic characterization of spaces goes back to Gel'fand and Kolmogorov \cite{KG}, who proved in 1939 that two compact topological spaces are homeomorphic if and only if the associative commutative algebras of continuous functions growing on them are isomorphic. Since then, similar characterizations have proved useful for def\/ining algebraically noncommutative spaces for instance.

A Lie algebra $\mathcal D(M)$ characterizes a manifold $M$ if for any manifold $N$, the Lie algebras~$\mathcal D(M)$ and~$\mathcal D(N)$ are isomorphic if and only if the manifolds $M$ and $N$ are dif\/feomorphic.
Lie-algebraic characterization of manifolds dates back to around 1950 with a result of L.E.~Pursell and M.E.~Shanks \cite{MR0064764} which states that the Lie-algebra of compactly supported vector f\/ields on a~manifold characterizes the latter manifold. Similar results have been obtained for various Lie subalgebras of ${\rm Vect}(M),$ the Lie algebra of vector f\/ields on $M$.

Rather than considering Lie subalgebras of vector f\/ields J.~Grabowski and N.~Poncin \cite{MR2027202,MR2172054,MR2346966} have characterized the manifold by means of the Lie algebra of linear dif\/ferential operators and the Poisson algebra of f\/iberwise polynomial functions on the cotangent bundle. Recently, they
extended, together with Kotov, their characterizations to the super context, see \cite{Grabowski-arXiv, MR2778234}.

Another way of generalizing Pursell--Shanks' result is the Lie algebraic characterization of vector bundles. Let us, for instance, point out a result of P.~Lecomte \cite{LEC} which shows that vector bundles are isomorphic if and only if their Lie algebras of inf\/initesimal automorphisms are isomorphic.

In the present paper we propose a hybrid approach in the sense that we give a characterization of vector bundles in terms of a Lie algebra which is not a Lie subalgebra of vector f\/ields. More precisely, we consider for a vector bundle $E\rightarrow M$, the Lie algebra $\mathcal D_{\mathcal E}(E)$ of \emph{homogeneous differential operators}, i.e., the Lie algebra generated by all dif\/ferential operators $D:{\rm C}^\infty(E)\to {\rm C}^\infty(E)$ which are eigenvectors with constant eigenvalues of $L_{{\mathcal E}}$, the Lie derivative in the direction of the Euler vector f\/ield.
We prove that the Lie algebra $\mathcal D_{\mathcal E}(E)$ and its Lie subalgebra $\mathcal D^1_{\mathcal E}(E)$ of f\/irst-order dif\/ferential operators both characterize the vector bundle $E\rightarrow M$.

This new Pursell--Shanks type result is thus remarkable in the sense that it is the whole f\/ibration that is characterized here, not only the base manifold or the total space of the vector bundle. However, since the proof relies on \cite[Theorem 1]{LEC}, our characterization inherits the same hypotheses and thus holds only for vector bundles of rank $n >1$ over base manifolds $M$ whose f\/irst space $\mathrm{H}^1(M,\mathbb{Z}/2)$ of the Cech cohomology vanishes.

In this work, we consider smooth manifolds which are assumed to be f\/inite-dimensional, Hausdorf\/f and second countable.

\section{The Lie algebra of homogeneous dif\/ferential operators}

Let $\pi: E\rightarrow M$ be a vector bundle. We denote by $\mathcal D^k(E,M)$ the space of $k$-th order dif\/ferential operators acting on ${\rm C}^\infty(E)$ and we def\/ine the set of \emph{homogeneous} dif\/ferential operators as
\[
\mathcal D_{\mathcal E}(E) =\bigcup_{k\geq0}\mathcal D^k_{\mathcal E}(E),
\]
where $\mathcal D^{k}_{\mathcal E}(E)$ is generated by $\{ T\in\mathcal D^k(E,M)\left|\right.\exists \,\lambda\in\mathbb{Z}: L_{\mathcal E} T=\lambda T\}$.
It follows from this def\/inition that $\mathcal D_{\mathcal E}$ is a $\mathbb{Z}$-graded space.
Moreover, $\mathcal D_{\mathcal E}$ is a Lie subalgebra of $\mathcal D(E,M)$. Actually, for any $T_1,T_2 \in \mathcal D_{_{\mathcal E}}(E)$ such that $L_{{\mathcal E}}(T_1) = \lambda T_1 $ and $L_{{\mathcal E}}(T_2)=\mu T_2$, one has
\[
L_{_{\mathcal E}}([T_1,T_2]) = (\lambda+\mu)[T_1,T_2].
\]
An important property of $\mathcal D_{\mathcal E}(E)$ is that the following relations hold:
\begin{gather*}
\mathcal D^k_{\mathcal E}(E)
  \subset   \mathcal D^{k+1}_{\mathcal E}(E) , \\
\mathcal D^k_{\mathcal E}(E) \mathcal D^l_{\mathcal E}(E)
  \subset   \mathcal D^{k+l}_{\mathcal E}(E)  , \\
{[\mathcal D^k_{\mathcal E}(E),\mathcal D^l_{\mathcal E}(E)]}
  \subset   \mathcal D^{k+l-1}_{\mathcal E}(E) .
\end{gather*}
Those relations show that, like $\mathcal D(E,M)$, the Lie subalgebra $\mathcal D_{\mathcal E}(E)$ is a \emph{quantum Poisson algebra}~\cite{MR2027202,MR2172054,MR2346966}. The \emph{basic algebra} of the quantum Poisson algebra $\mathcal D_{\mathcal E}(E)$ is the associative commutative algebra $\mathcal A(E)$ of f\/iberwise polynomial smooth functions on $E$, i.e.,
\[
\mathcal D^0_{\mathcal E}(E) \cong \mathcal A(E)=\oplus_{\lambda\in\mathbb{N}}\mathcal A^{\lambda}(E),
\]
where $\mathcal A^{\lambda}(E)=\{u\in{\rm C}^{\infty}(E):L_{_{\mathcal E}}u=\lambda u\}$.
Notice that $\mathcal A(E) \subset {\rm C}^{\infty}(E)$  is preserved by the elements of $D_{\mathcal E}(E)$. More precisely, if $u\in \mathcal A^\lambda(E)$ and $T \in \mathcal D_{\mathcal E}(E)$ is an eigenvector of $L_{_{\mathcal E}}$ corresponding to the eigenvalue $\mu,$ then
\[
T(u)\in\mathcal A^{\lambda+\mu}(E) .
\]
In the sequel, elements of $\mathcal D_{\mathcal E}(E)$ are usually considered as endomorphisms of the vector spa\-ce~$\mathcal A(E)$.

\subsection{A local description}

\begin{proposition}
A linear endomorphism $T : \mathcal A(E) \to \mathcal A(E)$ is an element of $\mathcal D_{\mathcal E}(E)$ if and only if in local coordinates $((x^i)_{1\leq i\leq m},(\xi_j)_{1\leq j\leq n})$ corresponding to a local trivialization of $E$, it reads
\[
T=T^{\alpha,\beta}\partial_{\alpha}\overline{\partial}_{\beta} ,
\]
where $\alpha$ and $\beta$ are multi-indices, $\partial_i=\frac{\partial}{\partial x^i}$, $\overline{\partial}_j=\frac{\partial}{\partial \xi_j}$ and all $T^{\alpha,\beta}$ are polynomials in $\xi_1, \dots, \xi_n$ with coefficients in ${\rm C}^{\infty}(M)$.
\end{proposition}
\begin{proof}
We will prove that if $L_{\mathcal E}T=\lambda T$ then every $T^{\alpha,\beta}$ is homogeneous of weight $\left|\beta\right|+\lambda$, i.e.,
\[
L_{\mathcal E}T^{\alpha\beta}= (\lambda+|\beta|) T .
\]
Observe that for all $f\in {\rm C}^{\infty}(E)$, one has
\begin{gather*}
(L_{\mathcal E}\overline{\partial}_i)f
  =   L_{\mathcal E}(\overline{\partial}_if)-\overline{\partial}_i(L_{\mathcal E}f)
  =   \sum_j\left(\xi^j\overline{\partial}_{ji}f\right) - \sum_j\left(\overline{\partial}_i(\xi^j\overline{\partial}_jf)\right) = - \overline{\partial}_if ,
\end{gather*}
so that
\begin{gather*}
L_{\mathcal E}(\overline{\partial}_{\beta})
  =   \sum_r\big(\overline{\partial}_{i_1}\cdots L_{\mathcal E}\overline{\partial}_{i_r}\cdots \overline{\partial}_{i_{\left|\beta\right|}}\big)
= -\left|\beta\right|\overline{\partial}_{\beta}.
\end{gather*}
Hence, for any $u\in\mathcal A(E)$ ,
\begin{gather*}
\lambda T^{\alpha,\beta}\partial_{\alpha}\overline{\partial}_{\beta}u
  =   L_{\mathcal E}\big(T^{\alpha,\beta}\partial_{\alpha}\overline{\partial}_{\beta}u\big) - T^{\alpha,\beta}\partial_{\alpha}\overline{\partial}_{\beta}(L_{\mathcal E}u)\\
\phantom{\lambda T^{\alpha,\beta}\partial_{\alpha}\overline{\partial}_{\beta}u}{}  =   (L_{\mathcal E}T^{\alpha,\beta})\partial_{\alpha}\overline{\partial}_{\beta}u + T^{\alpha,\beta}L_{\mathcal E}(\partial_{\alpha}\overline{\partial}_{\beta}u) - T^{\alpha,\beta}\partial_{\alpha}\overline{\partial}_{\beta}(L_{\mathcal E}u) \\
  \phantom{\lambda T^{\alpha,\beta}\partial_{\alpha}\overline{\partial}_{\beta}u}{} =   (L_{\mathcal E}T^{\alpha,\beta})\partial_{\alpha}\overline{\partial}_{\beta}u + T^{\alpha,\beta}L_{\mathcal E}(\partial_{\alpha}\overline{\partial}_{\beta})u
  =   L_{\mathcal E}T^{\alpha,\beta}\partial_{\alpha}\overline{\partial}_{\beta}u - \left|\beta\right|T^{\alpha,\beta}\partial_{\alpha}\overline{\partial}_{\beta}u .
\end{gather*}
It follows that
\[
\big(L_{\mathcal E}T^{\alpha,\beta}-\left(\left|\beta\right|+\lambda \right) T^{\alpha,\beta}\big)\partial_{\alpha}\overline{\partial}_{\beta} = 0,
\]
and the conclusion is immediate.
\end{proof}

\section[Properties of the quantum Poisson algebra ${\mathcal D}_{\mathcal E}(E)$]{Properties of the quantum Poisson algebra $\boldsymbol{{\mathcal D}_{\mathcal E}(E)}$}

\subsection{A non-singular quantum Poisson algebra}

Remember that the canonical isomorphism ${\mathcal D}_{\mathcal E}^0(E) \cong {\mathcal A}(E)$ is obtained by identifying a f\/iberwise polynomial smooth function $u$ with the operator of multiplication $\gamma_u : {\mathcal A}(E) \to {\mathcal A}(E) \text{, } v \mapsto uv$.

\begin{proposition}
The quantum Poisson algebra ${\mathcal D}_{\mathcal E}(E)$ is \emph{non-singular}, i.e.,
\[
\mathcal A(E) = \big[{\mathcal D}_{\mathcal E}^1(E) , \mathcal A(E)\big] .
\]
\end{proposition}

\begin{proof}
If $u\in\mathcal A^\lambda(E)$ for some $\lambda\neq0$, then
\[
\gamma_u = \big[\lambda^{-1}{\rm L}_{\mathcal E} , \gamma(u)\big] .
\]
If $u\in\mathcal A^0(E)$, then
\[
u = \sum_i [{\rm L}_{X_i^h}, \gamma_{\pi^*v_i}] ,
\]
if $u=\pi^*v$ and $v\in {\rm C}^{\infty}(M)$ reads
\[
v = \sum_i X_i \cdot v_i ,
\]
with $X_i\in {\rm Vect}(M)$ and $v_i\in {\rm C}^{\infty}(M)$. Since any $v \in {\rm C}^{\infty}(M)$ can be decomposed in this way and since~${\mathcal D}(M)$, the space of linear dif\/ferential operators of $M$, is non-singular (see \cite{MR2027202,MR2172054,MR2346966}), the conclusion follows.
\end{proof}

\subsection{The (quasi-)distinguishing property}

Denote by ${\mathcal S}_{\mathcal E}(E)$ the set $\{ \sigma(T) : T \in \mathcal D_{\mathcal E}(E) \}$ where $\sigma$ stands for the principal symbol operator. This set ${\mathcal S}_{\mathcal E}(E)$ is a subset in
\[
{\mathcal S}(E) = \{ \sigma(T) : T \in \mathcal D(E,M) \} \cong {\mathcal A}(T^*E) .
\] The Poisson bracket of ${\mathcal A}(T^*E)$ restricts to ${\mathcal S}_{\mathcal E}(E)$ and endows it with a commutative algebra structure.

The algebra $\mathcal D_{\mathcal E}(E)$ would be \emph{distinguishing} in the sense of \cite{MR2027202,MR2172054,MR2346966} if  for any $P \in {\mathcal S}_{\mathcal E}(E)$, we had
\begin{equation}\label{eq:distinguishing}
\left(\forall\, u \in\mathcal A(E), \exists\, n \in {\mathbb N} :
\underbrace {\{P,\{P,\dots\{P}_{n~\text{times}},u\}\}\} = 0 \right)
\ \Rightarrow \
\left( P \in {\mathcal S}_{\mathcal E}^0(E) = \mathcal A(E) \right).
\end{equation}
Obviously, this is not the case. Indeed, considering the trivial bundle $\rm{pr_1: \mathbb{R}^2 \to \mathbb{R}}$, it is easy to f\/ind a $P\in {\mathcal S}_{\mathcal E}(\mathbb{R}^2) \setminus {\mathcal S}_{\mathcal E}^0(\mathbb{R}^2)$ whose Hamiltonian vector f\/ield reads $H_P = \partial_{u^1}$ and for such a $P$, the right-hand side of condition (\ref{eq:distinguishing}) is trivially satisf\/ied since for any $\{P,-\}^n = H_P^n =  \partial_{u^1}^n$ while all elements in $\mathcal A_{\mathcal E}(E)$ are f\/inite order polynomial in $u_1$.

In \cite{MR2027202, MR2172054, MR2346966}, the fact that the Lie algebra of dif\/ferential operators on a manifold is non-singular and distinguishing is used to obtain that any isomorphism of such Lie algebras always have two nice properties : it is f\/iltered and its restriction to the basic algebras is an isomorphism of associative algebras. In order to obtain similar result for our Lie algebra, we need some ref\/inement here.

In general, see \cite{MR2027202, MR2172054, MR2346966}, for a quantum Poisson algebra ${\mathcal D}$ and its basic algebra ${\mathcal A}$ (whose unit is denoted by $1$), one refers to \emph{constants} for the elements of the image of the ground f\/ield~$K$ by the natural embedding $k \in K \to k1 \in {\mathcal A}$. The \emph{classical limit} of ${\mathcal D}$ is the Lie algebra $S({\mathcal D})=\oplus_{i\in\mathbb{Z}}{\mathcal D}^i/{\mathcal D}^{i-1}$ with the bracket given by
\[
\{ \sigma(D),\sigma(T)\}=
\left\{\begin{array}{ll}
\sigma([D,T]) , & \text{if} \ [D,T] \in {\mathcal D}^{\rm{ord}(D)+\rm{ord}(T)-1}\setminus{\mathcal D}^{\rm{ord}(D)+\rm{ord}(T)-2} ,\\
0 , &  \text{otherwise},
\end{array}\right.
\]
where the map $\sigma : {\mathcal D} \to S({\mathcal D})$ is the \emph{principal-symbol map} while $\rm{ord}(D) := i$ if $D \in {\mathcal D}^i\setminus{\mathcal D}^{i-1}$.
Note that  $\sigma({\mathcal A})={\mathcal A}=S_0({\mathcal D}).$
\begin{definition}
A quantum Poisson algebra ${\mathcal D}$ with basic algebra ${\mathcal A}$ and classical limit ${\mathcal S}$ is called \emph{quasi-distinguishing} if for any $P \in {\mathcal S}$,
\[
\{P,{\mathcal A}\}=0 \Rightarrow P \in {\mathcal A} .
\]
\end{definition}

We denote the \emph{centralizer} of ${\rm ad}_{\mathcal A}$ in ${\rm Hom}_K({\mathcal D},{\mathcal D})$ by
\[
{\mathcal C}({\mathcal D})= \{ \psi \in {\rm Hom}_K({\mathcal D},{\mathcal D}): \psi([D,u])=[\psi(D),u] , \; \forall\, D\in{\mathcal D}, \forall\, u\in{\mathcal A}\}.
\]

\begin{proposition}\label{prop:D ns and quasi}
If ${\mathcal D}$ is non-singular and quasi-distinguishing, then
\begin{enumerate}\itemsep=0pt
\item[$1)$]
for any $i \in {\mathbb N}$, $\{D\in {\mathcal D}: [D,{\mathcal A}]\subset {\mathcal D}^i\}= {\mathcal D}^{i+1}$ ;
\item[$2)$]
any $\psi \in {\mathcal C}({\mathcal D})$ preserves the filtration of ${\mathcal D}$ ;
\item[$3)$]
for any $\psi \in {\mathcal C}({\mathcal D})$ and any $u \in {\mathcal A}$, $\psi(u)=\psi(1)u$.
\end{enumerate}	
\end{proposition}

\begin{proof} The proof is similar to the corresponding one in \cite{MR2027202,MR2172054,MR2346966}. We give the main steps in order to show that it is enough to assume that ${\mathcal D}$ is quasi-distinguishing.

1.~If $D\in {\mathcal D}\setminus{\mathcal D}^{i+1}$ and $[D,{\mathcal A}]\subset {\mathcal D}^i$, then $\{\sigma(D), {\mathcal A}\}=0$. As ${\mathcal D}$ is quasi-distinguishing, $\sigma(D)\in {\mathcal A}$ and we have a contradiction since ${\mathcal A} \subset {\mathcal D}^{i+1}$.

2.~We proceed by induction. For any $D\in{\mathcal D}$, we have $\psi([D,{\mathcal A}])=[\psi(D),{\mathcal A}].$ In particular, when $D\in {\mathcal A}={\mathcal D}^0$, this equality reads $[\psi(D),{\mathcal A}]=0$ so that $\psi(D)\in{\mathcal A}.$
Now, if $\psi({\mathcal D}^i)\subset{\mathcal D}^i$ for some $i \in {\mathbb N}$, then for any $D\in{\mathcal D}^{i+1}$,
\[
[\psi(D),{\mathcal A}] = \psi([D,{\mathcal A}]) \in \psi({\mathcal D}^i) \subset {\mathcal D}^i.
\]
Using point 1) of the proposition, it follows that $\psi(D)\in{\mathcal D}^{i+1}.$

3.~First, note that  for any $D\in{\mathcal D}^1$ and any $u\in{\mathcal A}$, we have $[D,u^2]=2u[D,u]$. Applying $\psi$ to both sides, using point 2) of the proposition and the def\/ining property of ${\mathcal C}({\mathcal D})$, we get
\begin{gather*}
2\psi(u[D,u])
 =   \psi\big(\big[D,u^2\big]\big)
  =   \big[\psi(D),u^2\big]
  =   2u[\psi(D),u]
  =   2u\psi([D,u]).
\end{gather*}
In terms of the induced derivation $\hat{D}\in {\rm Der}({\mathcal A}):u\mapsto [D,u]$, we deduce that
\begin{equation}\label{eq:psi u}
\psi(u\hat{D}(u))=u\psi(\hat{D}(u)).
\end{equation}
In particular, for any $u,v,w \in {\mathcal A}$, we obtain
\[
\psi((u+w)v\hat{D}(u+w))=(u+w)\psi\big(v\hat{D}(u+w)\big).
\]
The above equality reduces to
\[
\psi\big(uv\hat{D}(w)\big)+\psi\big(wv\hat{D}(u)\big)=u\psi\big(v\hat{D}(w)\big)+w\psi\big(v\hat{D}(v)\big)
\]
and for $v=\hat{D}(w)$, it reads
\[
\psi\big(u\big(\hat{D}(w)\big)^2\big)+\psi\big(w\hat{D}(w)\hat{D}(u)\big)=u\psi\big(\big(\hat{D}(w)\big)^2\big)+w\psi\big(\hat{D}(w)\hat{D}(u)\big).
\]
As $\psi(w\hat{D}(w)\hat{D}(u)) = w\psi(\hat{D}(w)\hat{D}(u))$ in view of (\ref{eq:psi u}) with $D$ replaced by $\hat{D}(u)D$, we f\/inally obtain
\[
\psi\big(u\big(\hat{D}(w)\big)^2\big) = u\psi\big(\big(\hat{D}(w)\big)^2\big)
\]
for all $u,w \in {\mathcal A}$.
This shows that $[{\mathcal D}^1, {\mathcal A}]$ is a subset of the radical ${\rm Rad}({\mathcal J})$ of the ideal
\[
{\mathcal J}=\{ v\in{\mathcal A}: \psi(uv)=u\psi(v), \; \forall\, u\in{\mathcal A}\} \subset {\mathcal A}.
\]
Now since ${\mathcal D}$ is non-singular, we have $[{\mathcal D}^1, {\mathcal A}] = {\mathcal A}$ so that ${\rm Rad}({\mathcal J})$ and ${\mathcal J}$ both coincide with the whole algebra ${\mathcal A}$. The result follows immediately.
\end{proof}

\begin{definition}[\cite{MR2027202, MR2172054, MR2346966}]
A quantum Poisson algebra ${\mathcal D}$ is \textit{symplectic} if the constants are the only central elements.
\end{definition}

\begin{proposition}
Let ${\mathcal D_1}$ and ${\mathcal D_2}$ be non-singular, symplectic and quasi-distinguishing quantum Poisson algebras. Any isomorphism $\Phi :{\mathcal D_1}\to{\mathcal D_2}$ respects the filtration and the restriction ${\left.\Phi\right|}_{\mathcal A_1}: \mathcal A_1 \to \mathcal A_2$ is of the form
\[
{\left.\Phi\right|}_{\mathcal A_1}=k\Psi ,
\]
for some invertible element $k \in {\mathcal A_2}$, central in ${\mathcal D_2}$, and some isomorphism of $\mathbb{R}$-algebras $\Psi:\mathcal A_1 \to \mathcal A_2$.
\end{proposition}

\begin{proof}
For all $D\in {\mathcal D_2}$ and all $u\in {\mathcal A_1},$ $w\in {\mathcal A_2}$, we have
\begin{gather*}
\big(\Phi\circ\gamma_u\circ\Phi^{-1}\big)([D,w])
  =   \Phi\big[u\Phi^{-1}(D),\Phi^{-1}(w)\big]
  =   \big[\Phi\circ\gamma_u\circ\Phi^{-1}(D),w\big].
\end{gather*}
In view of point 2) and point 3) of Proposition~\ref{prop:D ns and quasi}, it follows that
\begin{equation}\label{formula:phi u phi-1}
\Phi\big(u\Phi^{-1}(w)\big)=\Phi\big(u\Phi^{-1}(1)\big)w\in {\mathcal A}_2
\end{equation}
so that $\Phi(u\Phi^{-1}(1))\in {\mathcal A}_2$.	As $\Phi^{-1}(1)$ is central in ${\mathcal D}_{1}$, it is a nonvanishing constant and thus $\Phi({\mathcal A_1})\subset{\mathcal A}_2$. Note that the same reasoning with $\Phi^{-1}$ yields $\Phi({\mathcal A_1})={\mathcal A}_2$.

Now assume that $\Phi({\mathcal D}_1^i)\subset {\mathcal D}_2^i$ for some $i \in {\mathbb N}$. For any $D\in{\mathcal D}_1^{i+1}$, we have
\[
[\Phi(D),{\mathcal A}_2]=\Phi([D,{\mathcal A}_1])\subset{\mathcal D}_2^{i} .
\]
Using point 1) of Proposition \ref{prop:D ns and quasi}, we obtain $\Phi(D)\in{\mathcal D}_2^{i+1}$ and thus $\Phi$ is f\/iltered.

Finally, setting $v=\Phi^{-1}(w)\in{\mathcal A}_1$ and $\lambda=\Phi^{-1}(1)$, equality (\ref{formula:phi u phi-1}) used twice yields
\begin{gather*}
\Phi(uv)
  =   \Phi(u\lambda)\Phi(v)
  =   \Phi(\lambda u)\Phi(v)
  =   \Phi\big(\lambda^2\big)\Phi(u)\Phi(v).
\end{gather*}
In particular, for $u=1$ and $v=\lambda,$ the relation above shows that the element $\Phi(\lambda^2)$, central in~${\mathcal D}_2$, is also invertible in~${\mathcal A}_2$. Setting $\kappa=\Phi(\lambda^2)^{-1}$, the association
\[
u\in {\mathcal A}_1 \ \mapsto \ \kappa^{-1}\Phi(u)\in {\mathcal A}_2,
\]
def\/ines an isomorphism of associative algebras and the conclusion is immediate.
\end{proof}

\begin{remark}
The same property holds for isomorphisms between the Lie subalgebras ${\mathcal D}_1^1$ and~${\mathcal D}_2^1$. The proof is exactly the same.	
\end{remark}

Notice that for a vector bundle $E \to M,$ the quantum Poisson algebra ${\mathcal D}_{\mathcal E}(E)$ is symplectic and quasi-distinguishing. These properties are easily shown from the local expression of the Poisson bracket in ${\mathcal S}_{\mathcal E}(E)$.

\begin{corollary}\label{cor:charact algb asso}
Let $\pi : E \to M$ and $\eta : F \to N$ be two vector bundles.
If $\Phi: {{\mathcal D}}_{\mathcal E}(E) \to {{\mathcal D}}_{\mathcal E}(F)$ is an isomorphism of
Lie algebras, then $\Phi$ maps ${\mathcal A}(E)$ on ${\mathcal A}(F)$ and the restriction ${\left.\Phi\right|}_{\mathcal A}(E): {\mathcal A}(E) \to {\mathcal A}(F)$ reads
\[
{\left.\Phi\right|}_{\mathcal A(E)}=\kappa\Psi ,
\]
for some $\kappa \in {\mathbb R}_0$ and some isomorphism of $\mathbb{R}$-algebras $\Psi:\mathcal A(E) \to \mathcal A(F)$.
\end{corollary}

\begin{remark}
The same property holds for isomorphisms between the Lie subalgebras ${\mathcal D}_{\mathcal E}^1(E)$ and~${\mathcal D}_{\mathcal E}^1(F)$.	
\end{remark}

\section{Isomorphisms induce graded isomorphisms}

\subsection{Isomorphisms are f\/iltered}

\begin{lemma}\label{lem:degree 0}
Let $\pi : E \to M$ and $\eta : F \to N$ be two vector bundles.
If $\Psi: {{\mathcal A}}(E) \to {{\mathcal A}}(F)$ is an isomorphism of
${\mathbb R}$-algebras, then
\[
\Psi\big({{\mathcal A}}^0(E)\big) = {{\mathcal A}}^0(F).
\]
\end{lemma}

\begin{proof}
For any nowhere vanishing $u \in {\mathcal A}^0(E)$, the function $u^{-1} : e \mapsto
\frac{1}{u(e)}$ also belongs to~${\mathcal A}^0(E)$. Since $\Psi$ is a homomorphism,
\[
{\mathcal A}^0(F) \ni 1_F = \Psi(1_E) = \Psi\big(u . u^{-1}\big) = \Psi(u) . \Psi\big(u^{-1}\big).
\]
This implies that $\Psi(u)$ and $\Psi(u^{-1})$ are both nonvanishing f\/iberwise
constant polynomials since their product is a nonvanishing f\/iberwise constant
polynomial.

Now, for any element $u \in {\mathcal A}^0(E)$, $u^2 + 1_{E} : e \mapsto u(e)^2 +
1$ is a nonvanishing element in ${\mathcal A}^0(E)$. It follows that
\[
{\mathcal A}^0(F) \ni \Psi\big(u^2+1_{E}\big) = \Psi(u).\Psi(u) + 1_{F},
\]
which shows that $\Psi(u)$ must be f\/iberwise constant.

We have proved $\Psi({{\mathcal A}}^0(E)) \subset {{\mathcal A}}^0(F)$. The conclusion
follows from applying the same arguments to the homomorphism~$\Psi^{-1}$.
\end{proof}

\begin{lemma}\label{lem:selection of low degree generators}
Let $\{u_1, u_2,\ldots,u_r\} \cup {{\mathcal A}}^0(E)$ be a system of generators of
the ${\mathbb R}$-algebra ${\mathcal A}(E)$. If $I \subset \{1,2,\ldots,r\}$ is the set of
indices $i$ for which $u_i$ is of order at most~$1$, then
\begin{equation}\label{set:generators refined}
\{u_i : i \in I\} \cup {{\mathcal A}}^0(E)
\end{equation}
is again a system of generators of ${\mathcal A}(E)$.
\end{lemma}

\begin{proof}\sloppy
It is obvious that the family (\ref{set:generators refined}) still generates
${\mathcal A}^0(E)$ and ${\mathcal A}^1(E)$. Since ${\mathcal A}^0(E)$ and~${\mathcal A}^1(E)$ generate the whole ${\mathcal A}(E)$, the lemma follows.
\end{proof}

\begin{proposition}\label{prop:iso are filtered}
Let $\pi : E \to M$ and $\eta : F \to N$ be two vector bundles.
Every isomorphism of ${\mathbb R}$-algebras $\Psi: {{\mathcal A}}(E) \to
{{\mathcal A}}(F)$ is filtered with respect to the filtrations of ${{\mathcal A}}(E)$ and
${{\mathcal A}}(F)$ associated with their gradings.
\end{proposition}

\begin{proof}
Since degrees $0$ and $1$ generate the whole ${\mathbb R}$-algebras, it is enough to
show
\begin{gather*}
\Psi\big({{\mathcal A}}^0(E)\big)   \subset   {{\mathcal A}}^0(F), \qquad
\Psi\big({{\mathcal A}}^1(E)\big)   \subset   {{\mathcal A}}^1(F) \oplus {{\mathcal A}}^0(F).
\end{gather*}
The f\/irst relation comes from Lemma~\ref{lem:degree 0}. Let us prove the
second one.

The ${{{\rm C}^\infty}(M)}$-module $\Gamma(E^*)$ of global sections of the dual bundle $\pi^*$ is (projective and) f\/initely generated. Let $\{\alpha_1,\alpha_2,\ldots,\alpha_r\}$ be a system of generators. Passing through the isomorphism $\Gamma(E^*) \cong {\mathcal A}^1(E)$ and adding f\/iberwise constant polynomials, we obtain a set of generators
\[
\{u_1, u_2, \ldots, u_r\} \cup {{\mathcal A}}^0(E)
\]
for the ${\mathbb R}$-algebra ${\mathcal A}(E)$. Since $\Psi$ is an isomorphism, the images
\[
\{\Psi(u_{1}),\Psi(u_{2}),\ldots,\Psi(u_{r})\} \cup {{\mathcal A}}^0(F)
\]
generate the whole ${\mathcal A}(F)$. In view of Lemma~\ref{lem:selection of low degree generators}, this set of generators can be ref\/ined into
\[
\{\Psi(u_{i}) : i \in I\} \cup {{\mathcal A}}^0(F),
\]
where $I \subset \{1,2,\ldots,r\}$ is the set of indices for which
$\Psi(u_{\alpha_i}) \in {\mathcal A}^1(F) \oplus {\mathcal A}^0(F)$. Since $\Psi^{-1}$ is
an isomorphism, every element in ${\mathcal A}^1(E)$ can now be decomposed as an
${\mathbb R}$-linear combination of the generating family
\[
\{u_{i} : i \in I\} \cup {{\mathcal A}}^0(E).
\]
The conclusion is immediate.
\end{proof}

\subsection{Isomorphisms induce graded isomorphisms}

For $k \in {\mathbb N}$, denote by ${\rm pr}_k$ the projection of f\/iberwise polynomial functions onto their homogeneous component of order $k$.

\begin{proposition}
If $\Psi: {{\mathcal A}}(E) \to {{\mathcal A}}(F)$ is an isomorphism of ${\mathbb R}$-algebras, then for any $k \in {\mathbb N}$, the map
\[
\widetilde{\Psi}_k  = {\rm pr}_k \circ {\left.\Psi\right|}_{{\mathcal A}^k(E)} : \ {\mathcal A}^k(E) \to {\mathcal A}^k(F)
\]
is a linear bijection.
\end{proposition}
\begin{proof}
First, $\widetilde{\Psi}_{k}$ is surjective. Given $w \in {\mathcal A}^k(F)$, there is an $u \in {\mathcal A}^k(E) \oplus \cdots \oplus  {\mathcal A}^1(E) \oplus {\mathcal A}^0(E)$ such that $w = \Psi(u)$. Using the fact  $\Psi$ is f\/iltered (in view of Proposition \ref{prop:iso are filtered}), it comes
\[
{\mathcal A}^k(F) \ni w = {\rm pr}_k \circ {\Psi(u)} = {\rm pr}_k \circ {\Psi({\rm pr}_k(u))} \in \widetilde{\Psi}_k\big({\mathcal A}^k(E)\big).
\]
Next, $\widetilde{\Psi}_{k}$ is also injective. If $u \in {\mathcal A}^k(E)$ satisf\/ies ${\rm pr}_k \circ \Psi(u) = 0$, then
$\Psi(u) \in {\mathcal A}^{k-1}(F) \oplus \cdots \oplus {\mathcal A}^{1}(F) \oplus {\mathcal A}^{0}(F)$. Since $\Psi^{-1}$ is f\/iltered (in view of Proposition \ref{prop:iso are filtered}), we get
\[
{\mathcal A}^{k}(E) \ni u = \Psi^{-1}(\Psi(u)) \in {\mathcal A}^{k-1}(E) \oplus \cdots \oplus {\mathcal A}^{1}(E) \oplus {\mathcal A}^{0}(E).
\]
Thus $u = 0$.
\end{proof}

\begin{corollary}\label{Psi iso => Psi tilde iso gradue}
If $\Psi: {{\mathcal A}}(E) \to {{\mathcal A}}(F)$ is an isomorphism of ${\mathbb R}$-algebras, then the map $\widetilde{\Psi} : {\mathcal A}(E) \to {\mathcal A}(F) $ whose restriction to each ${\mathcal A}^k(E)$ is given by $\widetilde{\Psi}_k$ is a graded isomorphism of ${\mathbb R}$-algebras.
\end{corollary}

\begin{proof}
The map $\tilde{\Psi}$ is obviously bijective, linear and graded. Moreover, it is a homomorphism of ${\mathbb R}$-algebras. Indeed, for any $u_1 \in {\mathcal A}^{k_1}(E)$ and $u_2 \in {\mathcal A}^{k_2}(E)$,
\begin{gather*}
\widetilde{\Psi}(u_1 u_2)
                =  {\rm pr}_{k_1+k_2}\circ \Psi(u_1 u_2)              ={\rm pr}_{k_1+k_2}(\Psi(u_1)\Psi(u_2)) \\
\phantom{\widetilde{\Psi}(u_1 u_2)}{} =  ({\rm pr}_{k_1}\circ \Psi(u_1)) ({\rm pr}_{k_2} \circ \Psi(u_2)) =  \widetilde{\Psi}(u_1) \widetilde{\Psi}(u_2),
\end{gather*}
where we have used the fact that the highest order terms of the product of f\/iberwise polynomial functions is the product of their highest order terms.
\end{proof}

\section{A Lie algebraic characterization}

After a technical lemma, we will be in a position to establish a new Pursell--Shanks type result for vector bundles, in terms of the Lie algebra generated by dif\/ferential operators that are eigenvalues of the Lie derivative in the direction of the Euler vector f\/ield associated with each vector bundle.

\begin{lemma}\label{lem1}
 If $0= j_{_a}^l u$, $u\in\mathcal{A}$, then
 \begin{gather}\label{(i)}
      u=\sum_{i=1}^N u_{i_0}\cdots u_{i_l}, 
\end{gather}
  where $u_{i_j}\in\mathcal{A}$ and $u_{i_j}(a)=0$,
                                $\forall\, (i,j)\in [1,N]\times[0,l]$.
\end{lemma}
\begin{proof}
Consider $u\in\mathcal{A}$ such that $j_{_{a}}^l(u)=0$ and $U$ a domain of trivialization  such that $\pi(a)\in U.$ One has over $U,$
\[
    \frac{d^{l+1}}{dt^{l+1}}u(a+th)
                 =\sum_{i_0,\dots,i_l}(\partial_{i_0\dots i_l}u)(a+th)h^{i_0}\cdots h^{i_l}.
\]
From this, it may be deduced by successive primitivations in $[0,t]$
  \[
    u(x)=\sum_{i_1,\dots,i_l}\big(x^{i_0}-a^{i_0}\big)\cdots\big(x^{i_l}-a^{i_l}\big)v(x),
  \]
where we put $x=a+th$.
In the right member of this equality, factors are local, but observe that one can write
  \[
    u=\big(1-\alpha^{l+1}\big)u+\alpha^{l+1}u,
  \]
with $\alpha=\pi^*(\alpha')$ where $\alpha'\in {\rm C}^\infty(M,\mathbb{R})$ equals to $1$ in a neighborhood $V\subset U$ containing $\pi(a)$ and with compact support in $U.$ The term $\alpha^{l+1}u$ becomes
   \[
     \sum_{i_1,\dots,i_l}\alpha\big(x^{i_0}-a^{i_0}\big)\cdots\alpha\big(x^{i_l}-a^{i_l}\big)v(x),
   \]
 and has the form \eqref{(i)}. For the other term, one proceeds in the same way by using the fact that
  \[
   1-\alpha^{l+1}=\big(1-\alpha^{l+1}\big)(1-\beta)^{l+1},
  \]
with $\beta=\pi^*(\beta')$ where $\beta'\in {\rm C}^\infty(M,\mathbb{R})$ is of compact support in a neighborhood $W\subset V$ containing $\pi(a),$ and in which it is equal to~$1$.
\end{proof}

\begin{corollary}\label{extention derivation}
  For all $D\in {\rm Der}(\mathcal A)$ there exists an unique differential operator
  $\widehat{D}\in \mathcal{D}^1(E,M)$ such that
    \[
      \widehat{D}(u)=D(u) \qquad \forall\, u\in\mathcal{A}.
    \]
\end{corollary}
\begin{proof}
We use a most general result stating that the corollary is true for a linear map
 \mbox{$D\!:\mathcal{A}\!\to\! \mathcal{A}$} which depends only on $1$-jets of its arguments (in the sense that $D$ vanishes on arguments of zero $1$-jet).
Indeed, this statement results from the fact that for any function $f\in {\rm C}^\infty(E,\mathbb{R})$, any point $a\in E $ and any integer $k\in\mathbb{N}$, there exists $u\in\mathcal{A}$ such that
   \[
     j_{_{a}}^k(f)=j_{_{a}}^k(u);
   \]
hence one can write
  \[
    \widehat{D}(f)(a)=D(u)(a).
  \]
Considering $D\in {\rm Der}(\mathcal{A})$ and $u\in\mathcal{A}$ such that $j_{_{a}}^1(u)=0,$
one has
  \[
    D(u)(a)=0.
  \]
Indeed, with the notations of Lemma~\ref{lem1}, one has
  \[
    D(u)= \sum_{i=1}^N u_{i_0}D(u_{i_1})+\sum_{i=1}^N u_{i_0}D(u_{i_1}).
  \]
Therefore $D(u)$ vanishes at $a$.
\end{proof}

\subsection{The main result}

The proof of the main result will rely on a theorem of \cite{LEC}, which we shall recall before stating out main result. As this theorem requires some hypotheses on both the base manifold and the rank of the vector bundles, our result will inherit those hypotheses.

\begin{theorem}[\cite{LEC}]\label{thm: lec}
Let $\pi: E \to M$ and $\eta: F \to N $ be two vector bundles of ranks $n,n' >1$. If the first space $\mathrm{H}^1(M,\mathbb{Z}/2)$ of the Cech cohomology of $M$ vanishes, then the two vector bundles are isomorphic if and only if the Lie algebras $\mathrm{Aut}(E)$ and $\mathrm{Aut}(F)$ of infinitesimal automorphisms are isomorphic.
\end{theorem}

\begin{theorem}\label{thm: main}
Under the same hypotheses as in Theorem {\rm \ref{thm: lec}}, two vector bundles $\pi: E \to M$ and $\eta: F \to N $ are isomorphic if and only if the Lie algebras ${\mathcal D}_{\mathcal E}(E)$ and ${\mathcal D}_{\mathcal E}(F)$
$($resp.\ ${\mathcal D}_{\mathcal E}^1(E)$ and~${\mathcal D}_{\mathcal E}^1(F))$
are isomorphic.
\end{theorem}

\begin{proof}
First, note that if $\Psi: {\mathcal A}(E) \to {\mathcal A}(F)$ is an isomorphism of graded ${\mathbb R}$-algebras ${\mathcal A}(E) \to {\mathcal A}(F)$, the induced isomorphism $\widehat{\Psi} : {\rm Der}(\mathcal A(E)) \to {\rm Der}(\mathcal A(E))$ transforms $0$-weight derivations of $\mathcal A(E)$ into $0$-weight derivations of $\mathcal A(F)$.
Then, the theorem is a consequence of Theorem~\ref{thm: lec}, if we prove that $\mathrm{Der}^0(\mathcal A(E))= {\rm Aut}(E)\left|_{\mathcal A(E)} \right.$.
On the one hand, the inclusion $\mathrm{Der}^0(\mathcal A(E))\supset {\rm Aut}(E)\left|_{\mathcal A(E)} \right.$ is obvious. On the other hand, for $D\in \mathrm{Der}^0(\mathcal A(E))$, Corollary \ref{extention derivation} yields $\widehat{D}\in {\rm Vect}(E)$ such that
$D=\widehat{D}\left|_{\mathcal A(E)} \right.$. As $[\mathcal E_E,\widehat{D}]\left|_{\mathcal A(E)} \right.=0$, it follows that $[\mathcal E_E,\widehat{D}]=0$.
\end{proof}

\begin{remark}
Along the way, we have also proved that {under the same hypotheses as above,} two vector bundles $\pi: E \to M$ and $\eta: F \to N $ are isomorphic if and only if the ${\mathbb R}$-algebras~${\mathcal S}_{\mathcal E}(E)$ and~${\mathcal S}_{\mathcal E}(F)$ (resp.\ ${\mathcal S}_{\mathcal E}^1(E)$ and ${\mathcal S}_{\mathcal E}^1(F)$) are isomorphic. Finally, note that without any hypo\-the\-sis, Corollary~\ref{cor:charact algb asso} and Lemma~\ref{lem:degree 0} ensure in view of Milnor's theorem that any isomorphism between~${\mathcal D}_{\mathcal E}(E)$ and~${\mathcal D}_{\mathcal E}(F)$ induces a dif\/feomorphism between the base manifolds $M$ and $N$.
\end{remark}

\subsection*{Acknowledgements}

We thank the referees for suggestions leading to improvements of the original paper.

\pdfbookmark[1]{References}{ref}
\LastPageEnding

\end{document}